\documentclass[12pt,a4paper]{article}
\usepackage[utf8]{inputenc}
\usepackage{amsmath}
\usepackage{amsfonts}
\usepackage{amssymb}
\usepackage{amsthm}
\usepackage{authblk}
\usepackage{inputenc}
\usepackage{longtable}

\allowdisplaybreaks

\begin{document}
\title{Subsets of colossally abundant numbers}
\author{Xiaolong Wu}
\affil{Ex. Institute of Mathematics, Chinese Academy of Sciences}
\affil{xwu622@comcast.net}
\date{March 11, 2019}
\maketitle

\begin{abstract}

    Let $G(n)=\sigma (n)/(n \log \log n )$. Robin made hypothesis that $G(n)<e^\gamma$ for all integer $n>5040$. This article divides all colossally abundant numbers in to three disjoint subsets CA1, CA2 and CA3, and shows that Robin hypothesis is true if and only if all CA2 numbers $>5040$ satisfy Robin inequality. 
\end{abstract}
\begin{center}\textbf{ \large Introduction}
\end{center}

Define $\rho(n)=\sigma(n)/n$, where $\sigma(n)=\sum_{d|n}d$ is the sum of divisor function. Define
\begin{equation}
G(n):=\frac{\rho(n)}{\log \log n}.
\end{equation}
Then Robin hypothesis is: all integers $n>5040$ satisfy Robin inequality
\begin{equation}\tag{RI}
G(n)<e^\gamma,
\end{equation}
where $\gamma$ is the Euler constant. Let
\begin{equation}
F(x,k):=\frac{\log (1+1/(x+x^2+\cdots+x^k ))}{\log x}.
\end{equation}
Define a set
\begin{equation}
E=\{{F(p,k) \, |\, prime\, p ,\,integer\,  k\geq 1}\}.
\end{equation}
Elements $\epsilon_i\in E$ are indexed in decreasing order. Elements in E are called critical parameters, For a given critical parameter $\epsilon_i$, we can construct a colossally abundant (abbreviate CA) number as follows: Define $x_k$ as the solution of
 \begin{equation}
 F(x_k,k)=\epsilon, \quad 1\leq k \leq K,
 \end{equation}
where K is the largest integer such that $x_K\leq 2$. For each prime define
\begin{equation}
a_p=\begin{cases}
k, & \text{if $x_{k+1}<p\leq x_k,\,$}\\
0, & \text{if $p>x_1$.}
\end{cases}
\end{equation}
and define
\begin{equation}
n_i:=\prod_p p^{a_p}.
\end{equation}
It can be proved that $n_i$ is a CA number, and $n_i$ will be called the CA number constructed from $\epsilon_i$. cf. [Broughan 2017] Section 6.3.
    For any integer $n\geq 2$, we will write $P(n)$ for the largest prime factor of n.
    
    We divide CA in to 3 disjoint subsets. Let $n_i$ be the CA number constructed from $\epsilon_i$, and p be the prime succeeding $P(n_i )$. 
    
    $n_i$ is called a CA1 number if $\log n_i<P(n_i )$. Theorem 1 shows $G(n_i )<G(n_{i-1} ),\forall \, n_i\in CA1,i\geq 3$.
    
    $n_i$ is called a CA2 number if $P(n_i )<\log n_i<p$. 
    
    $n_i$ is called a CA3 number if $p<\log n_i$. Let $n_j$ be the CA number constructed from F(p,1). Theorem 2 shows that $G(n_i )<G(n_j)$.
    
    Corollary 4 shows that Robin hypothesis is true if and only if all CA2 numbers $>5040$ satisfy (RI).

\begin{center}
             \textbf{Table 1. CA1 and CA2 numbers in the first 26 CA numbers}
\begin{longtable}{| r | r | r | r | r | r |}
\hline
\textbf{index i} & \textbf{$\log n_i$} & \textbf{$P(n_i)$} & is CA1? & is CA2? &$G(n_i)$\\
\hline
1	&0.6931	&2	&Y	&N	&-4.0926\\
\hline
2	&1.7918	&3	&Y	&N	&3.4294\\
\hline
3	&2.4849	&3	&Y	&N	&2.5634\\
\hline
4	&4.0943	&5	&Y	&N	&1.9864\\
\hline
5	&4.7875	&5	&Y	&N	&1.9157\\
\hline
6	&5.8861	&5	&N	&Y	&1.8335\\
\hline
7	&7.8320	&7	&N	&Y	&1.8046\\
\hline
8	&8.5252	&7	&N	&Y	&1.7910\\
\hline
9	&10.9231	&11	&Y	&N	&1.7512\\
\hline
10	&13.4880	&13	&N	&Y	&1.7331\\
\hline
11	&14.1812	&13	&N	&Y	&1.7277\\
\hline
12	&15.2798	&13	&N	&Y	&1.7235\\
\hline
13	&16.8892	&13	&N	&Y	&1.7179\\
\hline
14	&19.7224	&17	&N	&N	&1.7243\\
\hline
15	&22.6669	&19	&N	&Y	&1.7342\\
\hline
16	&25.8023	&23	&N	&Y	&1.7374\\
\hline
17	&26.4955	&23	&N	&Y	&1.7371\\
\hline
18	&29.8628	&29	&N	&Y	&1.7337\\
\hline
19	&33.2968	&31	&N	&Y	&1.7340\\
\hline
20	&35.2427	&31	&N	&Y	&1.7369\\
\hline
21	&36.3413	&31	&N	&Y	&1.7364\\
\hline
22	&39.9522	&37	&N	&Y	&1.7375\\
\hline
23	&43.6658	&41	&N	&N	&1.7380\\
\hline
24	&47.4270	&43	&N	&N	&1.7403\\
\hline
25	&48.1201	&43	&N	&N	&1.7406\\
\hline
26	&51.9703	&47	&N	&Y	&1.7430\\
\hline
\end{longtable}
\end{center}

    So, the smallest CA1 number is $n_1=2$; the smallest CA2 number is $n_6=360$; the smallest CA3 number is $n_{14}=367\, 567\, 200$.
        
    We next calculate the bounds of increment for $n_i\in CA3$. Let $p>10^8$ be the prime succeeding to $P(n_i )$. Assume $\epsilon_{i+1}=F(q,k)$ for some prime q and integer k. Then
    
    Theorem 3 shows a lower bound
     \[\frac{G(n_{i+1})}{G(n_i)} >\left(1-\frac{\log q}{3p^2 (\log p )^2}\right)^{-1}.\]
     
    Theorem 4 shows an upper bound
     \[\frac{G(n_{i+1})}{G(n_i)}<\exp\left(\frac{0.126\log q}{p(\log p)^3}\right).\]

    I checked the first 5 763 320 CA numbers (i.e. with the largest prime factor up to $10^8$). They contain 120 529 CA1 numbers, 5 565 CA2 numbers and 5 637 226 CA3 numbers.

\begin{center}
\textbf{ \large Main Content}
\end{center}

\noindent {\bfseries Lemma 1.}
\textit{ Let $\epsilon \in E$ be a critical parameter and $k\geq 1$ be an integer. Let $x_1$ and $x_k$ be defined by (4). Then 
\begin{equation}\tag{L1.1}
(x_k+\cdots +x_k^k)\log x_k\geq x_1 \log x_1+\left(1-\frac{1}{2x_1}\right)\left(\frac{\log x_1}{2}-\frac{\log x_k}{2}\right).
\end{equation}
\begin{equation}\tag{L1.2}
(x_k+\cdots +x_k^k)\log x_k< x_1 \log x_1+\frac{\log x_1}{2}-\frac{\log x_k}{2}+\frac{\log x_1}{4x_1}.
\end{equation}
(L1.1) and (L1.2) mean simple version:
\begin{equation}\tag{L1.1'}
(x_k+\cdots +x_k^k)\log x_k\geq x_1 \log x_1.
\end{equation}
\begin{equation}\tag{L1.2'}
(x_k+\cdots +x_k^k)\log x_k< x_1 \log x_1+\frac{\log x_1}{2}.
\end{equation}
}

\begin{proof}
By definition of $x_1$ and $x_k$, we have
\begin{equation}\tag{L1.3}
\frac{\log \left(1+\frac{1}{x_k+\cdots +x_k^k}\right)}{\log x_k}=\epsilon=\frac{\log \left(1+\frac{1}{x_1}\right)}{\log x_1}.
\end{equation}
\begin{equation}\tag{L1.4}
x_k^\epsilon =1+\frac{1}{x_k+\cdots +x_k^k},\quad x_1^\epsilon=1+\frac{1}{x_1}.
\end{equation}
Hence
\begin{align*}
\frac{x_k+\cdots +x_k^k}{x_1}&=\frac{x_1^\epsilon-1}{x_k^\epsilon-1}=\frac{e^{\epsilon \log x_1}-1}{e^{\epsilon \log x_k}-1}\\
&=\frac{\epsilon \log x_1+\frac{(\epsilon \log x_1)^2}{2!}+\cdots}{\epsilon \log x_k+\frac{(\epsilon \log x_k)^2}{2!}+\cdots}\\
&=\frac{\log x_1}{\log x_k}\left(\frac{1+\frac{\epsilon \log x_1}{2!}+\cdots}{1+\frac{\epsilon \log x_k}{2!}+\cdots}\right).\tag{L1.5}
\end{align*}
Compare
\begin{equation}\tag{L1.6}
\frac{1+\frac{\epsilon \log x_1}{2!}+\cdots}{1+\frac{\epsilon \log x_k}{2!}+\cdots}\, and \, 1+\epsilon\left(\frac{\log x_1}{2}-\frac{\log x_k}{2}+c\right),
\end{equation}
where c is a to-be-determined real parameter. 
\begin{align*}
H&:=\left(1+\frac{\epsilon \log x_1+\cdots}{2!}\right)\\
&\quad-\left(1+\frac{\epsilon \log x_k+\cdots}{2!}\right)\left(1+\epsilon\left(\frac{\log x_1}{2}-\frac{\log x_k}{2}+c\right)\right)\\
&=\sum_{j=1}^\infty\frac{(\epsilon \log x_1)^{j-1}}{j!}-\sum_{j=1}^\infty\frac{(\epsilon \log x_k)^{j-1}}{j!}\\
&\quad -\epsilon\left(\frac{\log x_1}{2}-\frac{\log x_k}{2}+c\right)\sum_{j=1}^\infty\frac{(\epsilon \log x_k)^{j-1}}{j!}\\
&=\sum_{j=2}^\infty\frac{(\epsilon \log x_1)^{j-1}}{j!}-\sum_{j=2}^\infty\frac{(\epsilon \log x_k)^{j-1}}{j!}\\
&\quad -\left(\frac{\log x_1}{2}-\frac{\log x_k}{2}+c\right)\sum_{j=2}^\infty\frac{\epsilon^{j-1} (\log x_k)^{j-2}}{(j-1)!}\\
&=-\epsilon c+\sum_{j=3}^\infty\frac{(\epsilon \log x_1)^{j-1}}{j!}\\
&\quad -\sum_{j=3}^\infty\frac{\epsilon^{j-1} (\log x_k)^{j-2}}{j!}\left(\log x_k+\frac{j \log x_1}{2}-\frac{j \log x_k}{2}+jc\right).\tag{L1.7}
\end{align*}
To prove (L1.1), set $c=0$. The lower bound of H is
\begin{align*}
H&=\sum_{j=3}^\infty \frac{\epsilon^{j-1}}{j!}\left( (\log x_1)^{j-1}-(\log x_k)^{j-2}\left(\log x_k+\frac{j(\log x_1}{2}-\frac{j \log x_k}{2}\right)\right)\\
&=\sum_{j=3}^\infty \frac{\epsilon^{j-1}}{j!}\left(((\log x_1)^{j-1}-(\log x_k)^{j-1})-\frac{j (\log x_k)^{j-2}}{2}(\log x_1-\log x_k)\right)\\
&=\sum_{j=3}^\infty \frac{\epsilon^{j-1}(\log x_1-\log x_k)}{j!}\left(\sum_{m=0}^{j-2}(\log x_1)^m(\log x_k)^{j-2-m}-\frac{j (\log x_k)^{j-2}}{2}\right)\\
&>0.\tag{L1.8}
\end{align*}
Combine (L1.5), (L1.7) and (L1.8), we have
\begin{equation}\notag
\frac{x_k+\cdots +x_k^k}{x_1}>\frac{\log x_1}{\log x_k}\left(1+\epsilon\left(\frac{\log x_1}{2}-\frac{\log x_k}{2}\right)\right).
\end{equation}
Since
\begin{equation}\notag
\epsilon=\frac{\log \left(1+\frac{1}{x_1}\right)}{\log x_1}>\frac{1}{x_1\log x_1}-\frac{1}{2x_1^2 \log x_1},
\end{equation}
we get
\begin{align*}
(x_k+\cdots +x_k^k)\log x_k&>x_1\log x_1\left(1+\epsilon\left(\frac{\log x_1}{2}-\frac{\log x_k}{2}\right)\right)\\
&>x_1\log x_1+\left(1-\frac{1}{2x_1}\left(\frac{\log x_1}{2}-\frac{\log x_k}{2}\right)\right). \tag{L1.9}
\end{align*}
That is, (L1.1) holds.\\
To prove (L1.2), we have from (L1.7)
\begin{equation}\tag{L1.10}
H<-\epsilon c+\sum_{j=3}^\infty \frac{(\epsilon \log x_1)^{j-1}}{j!}.
\end{equation}
The summation in (L1.10) can be simplified as
\begin{align*}
\sum_{j=3}^\infty \frac{(\epsilon \log x_1)^{j-1}}{j!}&=\sum_{j=2}^\infty \frac{(\epsilon \log x_1)^j}{(j+1)!}\\
&<\frac{(\epsilon \log x_1)^2}{6}\sum_{j=0}^\infty \frac{(\epsilon \log x_1)^j}{j!}=\frac{(\epsilon \log x_1)^2}{6}e^{\epsilon \log x_1}\\
&=\frac{(\epsilon \log x_1)^2}{6}x_1^\epsilon=\frac{(\epsilon \log x_1)^2}{6}\left(1+\frac{1}{x_1}\right).\tag{L1.11}
\end{align*}
By (L1.3), $\epsilon<1/(x_1 \log x_1)$, and we have
\begin{align*}
H&<-\epsilon c+\frac{(\epsilon \log x_1)^2}{6}\left(1+\frac{1}{x_1}\right)\\
&<\frac{\epsilon}{2}\left(-2c+\frac{\log x_1}{3x_1}\left(1+\frac{1}{x_1}\right)\right)\leq 0,\quad for\, c=\frac{\log x_1}{4x_1},x_1\geq 2.\tag{L1.12}
\end{align*}
Combine (L1.5), (L1.6) and (L1.12), we get
\begin{align*}
(x_k+\cdots +x_k^k)\log x_k&<x_1\log x_1\left(1+\epsilon\left(\frac{\log x_1}{2}-\frac{\log x_k}{2}+\frac{\log x_1}{4x_1}\right)\right)\\
&<x_1\log x_1+\frac{\log x_1}{2}-\frac{\log x_k}{2}+\frac{\log x_1}{4x_1}. \tag{L1.13}
\end{align*}
\end{proof}

\noindent {\bfseries Theorem 1.}
\textit{Let $i\geq 3$ be an integer and $n_i$ be a CA1 number, $p=P( n_i )$. Then 
\begin{equation}\tag{1.1}
G(n_i)<G(n_{i-1})\left(1-\left(\frac{\log q}{p\log p}\right)^2\right), \,if\, n_i/n_{i-1}=q.
\end{equation}
\begin{equation}\tag{1.2}
G(n_i)<G(n_{i-1})\left(1-\left(\frac{\log q}{p\log p}\right)^2\right)\left(1-\left(\frac{\log r}{p\log p}\right)^2\right), \,if\, n_i/n_{i-1}=qr.
\end{equation}
}
\begin{proof}
$n_i\in CA1$ means $\log n_i<p$.\\
1)  $n_i/n_{i-1}=q$. Assume $\epsilon_i=F(q,k)$ for some prime q and integer $k \geq 1$.
\begin{align*}
\frac{G(n_i)}{G(n_{i-1})}&=\frac{\rho(n_i)\log \log n_{i-1}}{\rho(n_{i-1})\log \log n_{i}}\\
&=\frac{\log \log n_i+\log\left(1-\frac{\log q}{\log n_i}\right)}{\log \log n_i}\left(1+\frac{1}{q+\cdots+q^k}\right)\\
&<\left(1-\frac{\log q}{\log n_i\log \log n_i}\right)\left(1+\frac{1}{q+\cdots+q^k}\right)\\
&<\left(1-\frac{\log q}{p\log p}\right)\left(1+\frac{1}{q+\cdots+q^k}\right).\tag{1.3}
\end{align*} 
By Lemma 1 (L1.1'), we have
\begin{equation}\tag{1.4}
(q+\cdots +q^k)\log q\geq x_1 \log x_1 \geq p \log p.
\end{equation}
Hence
\begin{equation}\tag{1.5}
\frac{G(n_i)}{G(n_{i-1})}<\left(1-\frac{\log q}{p\log p}\right)\left(1+\frac{\log q}{p\log p}\right)=1-\frac{(\log q)^2}{(p\log p)^2}.
\end{equation}
2) $n_i/n_{i-1}=qr$. Assume $\epsilon_i=F(q,k)=F(r,j)$ for some prime q, r and integer $k \geq 1, j\geq 1$. Then we have
\begin{align*}
&\frac{G(n_i)}{G(n_{i-1})}=\frac{\rho(n_i)\log \log n_{i-1}}{\rho(n_{i-1})\log \log n_{i}}\\
&=\frac{\log \log n_i+\log\left(1-\frac{\log q+\log r}{\log n_i}\right)}{\log \log n_i}\left(1+\frac{1}{q+\cdots+q^k}\right)\left(1+\frac{1}{r+\cdots+r^j}\right)\\
&<\left(1-\frac{\log q+\log r}{\log n_i\log \log n_i}\right)\left(1+\frac{1}{q+\cdots+q^k}\right)\left(1+\frac{1}{r+\cdots+r^j}\right)\\
&<\left(1-\frac{\log q}{p\log p}\right)\left(1-\frac{\log r}{p\log p}\right)\left(1+\frac{1}{q+\cdots+q^k}\right)\left(1+\frac{1}{r+\cdots+r^j}\right).\tag{1.6}
\end{align*}
By Lemma 1 (L1.1'), we have
\begin{equation}\tag{1.7}
(q+\cdots +q^k)\log q\geq p \log p,\quad (r+\cdots +r^j)\log r\geq p \log p.
\end{equation}
Hence we get
\begin{align*}
\frac{G(n_i)}{G(n_{i-1})}&<\left(1-\frac{\log q}{p\log p}\right)\left(1+\frac{\log q}{p\log p}\right)\left(1-\frac{\log r}{p\log p}\right)\left(1+\frac{\log r}{p\log p}\right)\\
&=\left(1-\frac{(\log q)^2}{(p\log p)^2}\right)\left(1-\frac{(\log r)^2}{(p\log p)^2}\right)\tag{1.8}
\end{align*}
\end{proof}

\noindent {\bfseries Corollary 1.}
\textit{Let $n_i> n_8=5040$ be a CA1 number. Let $n_j$ be the largest non-CA1 number below $n_i$. Then $G(n_i )<G(n_j )$.
}
\begin{proof}
The condition $n_i> n_8=5040$ guarantees the existence of $n_j$. By Theorem 1, we have
\begin{equation}
G(n_i )<G(n_{i-1} )<\cdots<G(n_{j+1} )<G(n_j ).
\end{equation}
\end{proof}

\noindent {\bfseries Corollary 2.}
\textit{Robin hypothesis is true if and only if all non-CA1 numbers $>5040$ satisfy (RI).
}
\begin{proof}
If one non-CA1 number $>5040$ fails (RI), then Robin hypothesis fails by definition. Conversely, if Robin hypothesis fails, then (RI) fails for a CA number $n_i>5040$, [NY 2014] Proposition 20. If $n_i\notin CA1$, then we are done. If $n_i\in CA1$, then by Corollary 1, there exists $n_j\notin CA1$, such that $G(n_i )<G(n_j )$. That is, (RI) fails for $n_j$.
\end{proof}

\noindent {\bfseries Lemma 2.}
\textit{Let $\epsilon\in E$ be a critical epsilon value. $x_k$ are solutions of  
\begin{equation}\tag{L2.1}
F(x_k,k)=\epsilon.\quad  k\geq 1.
\end{equation}
Then $g(t)=g_\epsilon (t):=t^\epsilon/\log \log t$  has a unique minimum, say $t_0$, and $t_0$ satisfies
\begin{equation}\tag{L2.2}
x_1+\frac{1}{2}-\frac{1}{2 \log x_1}<\log t_0<x_1+\frac{1}{2}-\frac{1}{12x_1}+\frac{1}{24x_1^2},\quad   \forall\, x_1\geq 2.
\end{equation}
}
\begin{proof}
Take derivative,
\begin{align*}
g'(t)&=\frac{\epsilon t^{\epsilon-1}\log \log t-t^\epsilon \frac{1}{t \log t}}{(\log \log t)^2}\\
&=\frac{t^{\epsilon-1}}{\log t (\log \log t)^2}(\epsilon \log t \log \log t-1).\tag{L2.3}
\end{align*}
Define 
\begin{equation}\tag{L2.4}
f(t):=\epsilon \log t \log \log t-1.
\end{equation}
It is obvious that f(t) monotonically increases for $t\in (e,\infty)$, negative near e and positive when t sufficiently large. So $f(t)$ has a unique zero $t_0$. $g(t)$ attains minimum at $t_0$. Note $x_1$ is the solution of $F(x_1,1)=\log (1+1/x_1 )/\log x_1=\epsilon$, Write $t=x_1+1/2+d$, where $d=-1/(2 \log x_1)$. We have
\begin{align*}
f\left( e^{x_1+\frac{1}{2}+d}\right)&=\frac{\log (1+1/x_1)}{\log x_1}\left(x_1+\frac{1}{2}+d\right)\log \left(x_1+\frac{1}{2}+d\right)-1\\
&=\frac{\log (1+\frac{1}{x_1})}{\log x_1}\left(x_1+\frac{1}{2}+d\right)\left(\log x_1+\log\left(1+\frac{\frac{1}{2}+d}{x_1} \right)\right)-1\\
&<\left(\frac{1}{x_1}-\frac{1}{2x_1^2}+\frac{1}{3x_1^3}\right)\left(x_1+\frac{1}{2}+d\right)\\
&\quad \times \left(1+\frac{1}{2x_1 \log x_1}-\frac{1}{2x_1(\log x_1)^2}\right)-1\\
&=\left(1-\frac{1}{2x_1\log x_1}+\frac{1}{12x_1^2}+\frac{1}{4x_1^2\log x_1}+\frac{1}{6x_1^3}-\frac{1}{6x_1^3\log x_1}\right)\\
&\quad \times \left(1+\frac{1}{2x_1\log x_1}-\frac{1}{2x_1(\log x_1)^2}\right)-1\\
&<0,\quad \forall\,x_1\geq 2. \tag{L2.5}
\end{align*}
So we get the left inequality of (L2.2). For the right inequality, we have
\begin{align*}
f\left(e^{1/\log (1+1/x_1)}\right)&=\frac{\log (1+1/x_1)}{\log x_1}\frac{1}{\log (1+1/x_1)}\log \left(\frac{1}{\log(1+1/x_1)}\right)-1\\
&=\frac{1}{\log x_1}\log \left(x_1+\frac{1}{2}-\frac{1}{12x_1}+\frac{1}{24x_1^2}-\cdots\right)-1>0,\tag{L2.6}
\end{align*}
here the expansion of $(\log (1+1/x_1 ) )^{-1}= x+\frac{1}{2}-\frac{1}{12x_1}+\frac{1}{24x_1^2}-\cdots$ is calculated term wise from the formula $\log \left(1+\frac{1}{x}\right)=\frac{1}{x}-\frac{1}{2x^2}+\frac{1}{3x^3}-\cdots$. So we have
\begin{align*}
\log t_0&<\frac{1}{\log (1+\frac{1}{x_1})}=x_1+\frac{1}{2}-\frac{1}{12x_1}+\frac{1}{24x_1^2}-\cdots\\
&<x_1+\frac{1}{2}-\frac{1}{12x_1}+\frac{1}{24x_1^2}.\tag{L2.7}
\end{align*}
\end{proof}

\noindent {\bfseries Lemma 3.}
\textit{ Let $\epsilon\in E$ be a critical epsilon value. Let u and $u_1<u_2$ be positive reals. Then $h(u):= e^{\epsilon u}/\log u$  has a unique minimum at $u_0=u_0 (\epsilon)$ implicitly defined by $\epsilon=1/(u_0 \log u_0)$. Assume $u_0>40$. Write
\begin{equation}\tag{L3.1}
h_0:=h(u_0)=\frac{e^{1/\log u_0}}{\log u_0}.
\end{equation} 
1) For $u_0-\frac{1}{2}<u_1<u_0$,
\begin{equation}\tag{L3.2}
\frac{h(u_2)}{h_0} <1+0.2532\frac{u_0-u_1}{u_0^2\log u_0}+0.5162\frac{(u_0-u_1)^2}{u_0^2(\log u_0)^2}.
\end{equation}
2) For $u_0<u_2<u_0\log u_0$,
\begin{equation}\tag{L3.3}
\frac{h(u_2)}{h_0}>1+\frac{(u_2-u_0)^2}{2u_0^2\log u_0}-\frac{(u_2-u_0)^2}{2u_0^2(\log u_0)^2}.
\end{equation}
3) For $u_0<u_1<u_2,\,u_2-u_1<\log u_0$,
\begin{equation}\tag{L3.4}
\frac{h(u_2)}{h(u_1)} >1+0.3337\frac{(u_2-u_1)^2}{u_0^2 (\log u_0)^2}.
\end{equation}
}
\begin{proof}
We have
\begin{align*}
h(u)&=\frac{e^{\epsilon u}}{\log u}=\frac{e^{\epsilon u_0}e^{\epsilon (u-u_0)}}{\log u_0+\log (u/u_0)}=\frac{e^{\epsilon u_0}}{\log u_0}\cdot\frac{e^{\epsilon (u-u_0)}}{1+\frac{\log (u/u_0)}{\log u_0}}\\
&=h_0\left(\sum_{i=0}^\infty\frac{(\epsilon(u-u_0))^i}{i!}\right)\left(\sum_{j=0}^\infty\left(\frac{-\log (u/u_0)}{\log u_0}\right)^j\right).\tag{L3.5}
\end{align*}
1) When $u_0-\frac{1}{2}<u_1<u_0$, we have $\log (u_1/u_0)<0$. Hence
\begin{align*}
\frac{-\log (u_1/u_0)}{\log u_0}&=\frac{-\log(1-(u_0-u_1)/u_0)}{\log u_0}\\
&<\frac{u_0-u_1}{u_0 \log u_0}\sum_{k=1}^\infty\frac{1}{k}\left(\frac{u_0-u_1}{u_0}\right)^{k-1}\\
&<\frac{u_0-u_1}{u_0\log u_0}\left(1+\frac{1}{4u_0}\sum_{k=2}^\infty \left( \frac{1}{2u_0}\right)^{k-2}\right)\\
&<\frac{u_0-u_1}{u_0\log u_0}\left(1+\frac{0.2532}{u_0}\right).\tag{L3.6}
\end{align*}
\begin{align*}
\frac{h(u_1)}{h_0} &= \left(\sum_{i=0}^\infty\frac{(-\epsilon(u_0-u_1))^i}{i!}\right)\left(\sum_{j=0}^\infty \left(\frac{u_0-u_1}{u_0\log u_0}\left(1+\frac{0.2532}{u_0}\right)\right)^j\right)\\
&<\left(1-\frac{u_0-u_1}{u_0\log u_0}+\frac{(u_0-u_1)^2}{2u_0^2(\log u_0)^2}\right)\\
&\quad \times\left(1+\frac{u_0-u_1}{u_0\log u_0}+0.2532\frac{u_0-u_1}{u_0^2\log u_0}+1.0162\frac{(u_0-u_1)^2}{u_0^2(\log u_0)^2}\right)\\
&<1+0.2532\frac{u_0-u_1}{u_0^2\log u_0}+0.5162\frac{(u_0-u_1)^2}{u_0^2(\log u_0)^2}.\tag{L3.7}
\end{align*}
2) When $u_0<u_2<u_0\log u_0$, we have $\log (u_2/u_0 )>0$.
\begin{equation}\tag{L3.8}
\frac{-\log(u_2/u_0)}{\log u_0}=\frac{-\log(1+(u_2-u_0)/u_0)}{\log u_0}>-\frac{u_2-u_0}{u_0\log u_0}+\frac{(u_0-u_1)^2}{2u_0^2\log u_0}
\end{equation}
\begin{align*}
\frac{h(u_2)}{h_0} &> \left(\sum_{i=0}^\infty\frac{(\epsilon(u_2-u_0))^i}{i!}\right)\left(\sum_{j=0}^\infty \left(-\frac{u_2-u_0}{u_0\log u_0}+\frac{(u_2-u_0)^2}{2u_0^2\log u_0}\right)^j\right)\\
&>\left(\sum_{i=0}^\infty \frac{(u_2-u_0)^i}{i!u_0^i(\log u_0)^i}\right)\left(1-\frac{u_2-u_0}{u_0\log u_0}+\frac{(u_2-u_0)^2}{2u_0^2(\log u_0)^2}\right)\\
&>\left(1+\frac{u_2-u_0}{u_0\log u_0}+\frac{(u_2-u_0)^2}{2u_0^2(\log u_0)^2}\right)\left(1-\frac{u_2-u_0}{u_0\log u_0}+\frac{(u_2-u_0)^2}{2u_0^2(\log u_0)^2}\right)\\
&>1+\frac{(u_2-u_0)^2}{2u_0^2\log u_0}-\frac{(u_2-u_0)^2}{2u_0^2(\log u_0)^2}.\tag{L3.9}
\end{align*}
3) Write $u_2=u_1+a$ for some real $a<\log u_0$.
\begin{align*}
\frac{h(u_1)}{h(u_1)} &=e^{\epsilon (u_2-u_1)}\frac{\log u_1}{\log u_2}=e^{a/(u_0\log u_0)}\frac{\log u_1}{\log (u_1+a)}\\
&>e^{a/(u_0\log u_0)}\frac{\log u_0}{\log (u_0+a)}.\tag{L3.10}
\end{align*}
Since
\begin{align*}
\frac{\log u_0}{\log (u_0+a)}&=\frac{\log u_0}{\log u_0+\log (1+a/u_0)}=\frac{1}{1+\log (1+a/u_0)/ \log u_0}\\
&=\sum_{i=0}^\infty \left(\frac{-\log(1+a/u_0)}{\log u_0}\right)^i\\
&>1-\frac{\log (1+a/u_0)}{\log u_0}=1-\frac{1}{\log u_0}\sum_{i=1}^\infty\frac{(-1)^{i-1}}{i}\left(\frac{a}{u_0}\right)^i\\
&>1-\frac{a}{u_0\log u_0}+\frac{a^2}{2u_0^2\log u_0}-\frac{a^3}{3u_0^3\log u_0}.\tag{L3.13}
\end{align*}
and
\begin{equation}\tag{L3.14}
e^{a/(u_0\log u_0)}=\sum_{i=0}^\infty \frac{1}{i!}\left(\frac{a}{u_0\log u_0}\right)^i>1+\frac{a}{u_0\log u_0}+\frac{a^2}{2u_0^2\log u_0},
\end{equation}
we have
\begin{align*}
\frac{h(u_2)}{h(u_1)} &>\left(1+\frac{a}{u_0\log u_0}+\frac{a^2}{2u_0^2\log u_0}\right)\\
&\quad \times\left(1-\frac{a}{u_0\log u_0}+\frac{a^2}{2u_0^2\log u_0}-\frac{a^3}{3u_0^3\log u_0}\right)\\
&>1+\frac{a^2}{u_0^2\log u_0}\left(\frac{1}{2}-\frac{1}{2\log u_0}-\frac{a}{3u_0}\right)>1+\frac{0.3337a^2}{u_0^2\log u_0}.\tag{L3.15}
\end{align*}
\end{proof}

\noindent {\bfseries Lemma 4.}
\textit{Assume $g(t)=t^\epsilon/\log \log t$ takes minimum at  $t_0=t_0 (\epsilon)$. Assume $\log t_0>40$. Let N and $N_1$ be positive integers.\\
1) If $\log t_0-\frac{1}{2}<\log N<\log t_0$ and $\log t_0+2<\log N_1$, then 
\begin{equation}\tag{L4.1}
g(N_1)>g(N)\left(1+\frac{2.754}{(\log t_0)^2\log \log t_0}\right).
\end{equation}
2) If $\log t_0<\log N<\log N_1$ and $\log N_1-\log N<\log \log t_0$, then
\begin{equation}\tag{L4.2}
g(N_1)>g(N)\left(1+\frac{0.3337(\log N_1-\log N)^2}{(\log t_0)^2\log \log t_0}\right).
\end{equation}
}
\begin{proof}
Write $u=\log t,\,u_0=\log t_0$, $h(u)=g(t)$, $h_0=\log t_0$. By Lemma 3 (L3.1) and (L3.2), we have
\begin{align*}
\frac{g(N_1)-g(N)}{h_0}&>\frac{(\log N_1-u_0)^2}{u_0^2\log u_0}-\frac{(\log N_1-u_0)^2}{2u_0^2(\log u_0)^2}\\
&\quad -0.2532\frac{u_0-\log N}{u_0^2\log u_0}-0.5162\frac{(u_0-\log N)^2}{u_0^2(\log u_0)^2}\\
&>\frac{4-0.1266}{u_0^2\log u_0}-\frac{2+0.1291}{u_0^2(\log u_0)^2}>\frac{3.2962}{u_0^2\log u_0}.\tag{L4.3}
\end{align*}
By Lemma 3 (L3.1)
\begin{equation}\tag{L4.4}
\frac{g(N)}{h_0}<1+\frac{0.2532\times 0.5}{40^2\log 40}+\frac{0.5162\times 0.5^2}{40^2(\log 40)^2}=1.00002738,
\end{equation}
we have
\begin{equation}\tag{L4.5}
g(N_1)>g(N)\left(1+\frac{3.2962\,h_0}{g(N)u_0^2\log u_0}\right)>g(N)\left(1+\frac{3.2961}{u_0^2\log u_0}\right).
\end{equation}
2) follows from Lemma 3 (L3.3).
\end{proof}

\noindent {\bfseries Theorem 2.}
\textit{Let $n_i$ be CA3. Let p be the prime succeeding $P(n_i )$, $n_j$ be the CA number constructed from $\epsilon_j=F(p,1)$. then
\begin{equation}\tag{2.1}
G(n_j )>G(n_i )\left(1+\frac{3.2961}{(\log t_0)^2 \log \log t_0}\right),
\end{equation} 
where $t_0$ is defined as in Lemma 4.
}
\begin{proof}
$n_i\in CA3$ means $p<\log n_i$. By definition of CA numbers, we have
\begin{equation}\notag
\frac{\rho(n_i)}{n_i^{\epsilon_j}}\leq \frac{\rho(n_j)}{n_j^{\epsilon_j}}. 
\end{equation}
\begin{equation}\tag{2.2}
\frac{G(n_j)}{G(n_i)}=\frac{\rho(n_j) n_j^{\epsilon_j}}{n_j^{\epsilon_j}\log \log n_j}\cdot \frac{n_i^{\epsilon_j}\log \log n_i}{\rho(n_i) n_i^{\epsilon_j}}\geq\frac{n_j^{\epsilon_j}\log \log n_i}{n_i^{\epsilon_j}\log \log n_j}=\frac{g(n_j)}{g(n_i)},
\end{equation}
where $g(t)=t^{\epsilon_j}/\log \log t$ . By Lemma 2,  $g(t)$ attains minimum at $t_0$, and
\begin{equation}\tag{2.3}
p+\frac{1}{2}-\frac{1}{2\log p}<\log t_0<p+\frac{1}{2}-\frac{1}{12p}+\frac{1}{24p^2}.
\end{equation}
The smallest CA3 number $n_{14}$ can be directly checked. So, we may start from the next CA3 number $n_{23}$. That is, we may assume $n_i\geq n_{23}$ with $p\geq 43$ and $\log p \geq 3.76$.\\
Case 1) $\log n_i< \log t_0$. In this case, we have $\log t_0<p+\frac{1}{2}$ by (2.3). Hence $\log t_0-\frac{1}{2}<p<\log n_i <\log t_0$, and 
\begin{equation}\tag{2.4}
 \log n_j-\log t_0>\log n_j-\log n_i-\frac{1}{2}>\log p-\frac{1}{2}>2.
\end{equation}
Hence the conditions of Lemma 4 (L4.1) are satisfied and (2.1) holds.\\
Case 2) $\log n_i \geq \log t_0$. In this case, $\log n_j-\log n_i<\log t_0$. So by Lemma 4 (L4.2), we have
\begin{equation}\tag{2.5}
G(n_j)>G(n_i)\left(1+\frac{0.3337(\log n_j-\log n_i)^2}{(\log t_0)^2\log \log t_0}\right).
\end{equation}
Since $0.3337(\log n_j-\log n_i)^2\geq 0.3337\times (\log 43)^2>3.2961$, (2.1) holds.
\end{proof}

\noindent {\bfseries Corollary 3.}
\textit{Let $n_i$ be a CA3 number. Then there exists $n_j\in CA2$ such that $n_i<n_j$. If  $n_j$ is the smallest CA2 number above $n_i$, then $G(n_i )<G(n_j )$.
}
\begin{proof}
There are infinite CA1 numbers n, i.e. $\log n<P(n)$, [CNS 2012] Theorem 7. Let $n_k$ be the smallest such number above $n_i$. 

    We claim that $n_{k-1}$ is CA2. $n_{k-1}$ is not CA1 by minimality of $n_k$. If $n_{k-1}$ were CA3, there would exist a prime p such that $P(n_{k-1} )<p<\log n_{k-1}$. Then we would have 
\begin{equation}\tag{C3.1}
    \log n_k=\log n_{k-1}+\log (n_k/n_{k-1} )>\log n_{k-1} >p\geq P(n_k ).
\end{equation}    
This contradicts to $n_k\in CA1$. So $n_{k-1}\in CA2$ and we proved the existence of $n_j$. 

    Write $p_r=P(n_i )$, $p_s=P(n_j )$. Let $n_{i_m}$ be the CA number generated from parameter $F(p_m,1)$, $r<m\leq s$. Since $n_k$ is the smallest CA1 number above $n_i$, and $n_j<n_k$ is the smallest CA2 number above $n_i$, all $n_{i_m}<n_j$ are CA3. By Theorem 2, we have
\begin{equation}\tag{C3.2}
G(n_i )<G(n_{i_{r+1}} )<\cdots<G(n_{i_s})=G(n_j ).
\end{equation}
\end{proof}

\noindent {\bfseries Corollary 4.}
\textit{Robin hypothesis is true if and only if all CA2 numbers $>5040$ satisfy (RI).
}
\begin{proof}
If one CA2 number $>5040$ fails (RI), then Robin hypothesis fails by definition. Conversely, if (RI) fails, then by Corollary 2, (RI) fails for a non-CA1 number $n_i>5040$. If $n_i\in CA2$, then we are done. If $n_i\notin CA2$, then by Corollary 3, there exists $n_j\in CA2$, such that $G(n_i )<G(n_j )$. That is, (RI) fails for $n_j$.
\end{proof}

 Under assumption of Theorem 2, is $G(n_i )<G(n_{i+1} )$? Let $\epsilon_{i+1}=F(q,k)$. If $q \geq 3$, Theorem 3 proves $G(n_i )<G(n_{i+1} )$. The case $q=2$ is open. Theorem 3 also shows a lower bound for $G(n_{i+1} )/G(n_i )$ .

\noindent {\bfseries Theorem 3.}
\textit{Let $n_i$ be CA3. Let p be the prime succeeding $P(n_i )$. \\ 
1) If $\epsilon_{i+1}=F(q,k)$, $q\geq 3$, then $n_{i+1}/n_i=q$, and $G(n_i)<G(n_{i+1})$.\\
2) If $\epsilon_{i+1}=F(q,k)$, $q\geq 23$, then $n_{i+1}/n_i=q$, and
\begin{equation}\tag{3.1}
G(n_i)<G(n_{i+1})\left(1-\frac{(\log q)^2}{3p^2\log p}\right). 
\end{equation}
3) If $\epsilon_{i+1}=F(q,k)=F(r,j)$ , $q\geq 23,r\geq 23$, then $n_{i+1}/n_i=qr$, and
\begin{equation}\tag{3.2}
G(n_i)<G(n_{i+1})\left(1-\frac{(\log q)^2}{3p^2\log p}\right)\left(1-\frac{(\log r)^2}{3p^2\log p}\right). 
\end{equation}
}
\begin{proof}
I numerically checked for all CA3 numbers $n_i$ with $i<10\, 000$. $1)-3)$ all hold. So we may assume $i \geq 10\, 000$, and hence $p>103\, 049$. Since  $n_i\in CA3$, we have $p<\log n_i$.\\
1) and 2). Since $\epsilon_{i+1}=F(q,k)$, we have $n_{i+1}=n_i q$. Compare $G(n_i )$ and $G(n_{i+1} )$, we have 
\begin{align*}
\frac{G(n_i )}{G(n_{i+1} )}&=\frac{\rho(n_i)\log \log n_{i+1}}{\rho(n_{i+1})\log \log n_i}\\
&=\frac{\log(\log n_i+\log q)}{\log \log n_i}\left(\frac{q+\cdots+q^k}{1+q+\cdots+q^k}\right)\\
&=\frac{\log \log n_i+\log \left(1+\frac{\log q}{\log n_i}\right)}{\log \log n_i}\left(1-\frac{1}{1+q+\cdots+q^k}\right).\tag{3.3}
\end{align*}
Since $\log n_i>p>x_1$, where $x_1$ is defined by (4), Lemma 1 (L1.2’) means
\begin{align*}
\frac{G(n_i )}{G(n_{i+1} )}&<\left(1+\frac{\log \left(1+\frac{\log q}{p}\right)}{\log p}\right)\left(1-\frac{1}{1+\frac{1}{\log q}\left(p \log p+\frac{\log p}{2}\right)}\right)\\
&<\left(1+\frac{(\log q) \left(1-\frac{\log q}{2p}+\frac{(\log q)^2}{3p^2}\right)}{p\log p}\right)\left(1-\frac{1}{1+\frac{p \log p}{\log q}+\frac{\log p}{2\log q}}\right)\\
&=1+\frac{(p \log p+\frac{\log p}{2}) \left(1-\frac{\log q}{2p}+\frac{(\log q)^2}{3p^2}\right)-p \log p}{p\log p \left(1+\frac{p \log p}{\log q}+\frac{\log p}{2\log q}\right)}\\
&<1+\frac{\left(p+\frac{1}{2}\right)\left(1-\frac{\log q}{2p}+\frac{(\log q)^2}{3p^2}\right)-p}{p \left(1+\frac{p \log p}{\log q}+\frac{\log p}{2\log q}\right)}.\tag{3.4}
\end{align*}
\begin{align*}
&\left(p+\frac{1}{2}\right)\left(1-\frac{\log q}{2p}+\frac{(\log q)^2}{3p^2}\right)-p\\
&=\frac{1}{2}-\frac{\log q}{2}+\frac{(\log q)^2}{3p}-\frac{\log q}{2p}\left(\frac{1}{2}-\frac{\log q}{3p}\right)\\
&<\frac{1}{2}-\frac{\log q}{2}+\frac{(\log q)^2}{3p}. \tag{3.5}
\end{align*}
When $q\geq 3$, the expression in (3.5) is negative, so (3.4) means $G(n_i )<G(n_{i+1} )$. That is, 1) is true. Now for 2) we have
\begin{align*}
\frac{G(n_i )}{G(n_{i+1} )}&<1+\frac{\frac{1}{2}-\frac{\log q}{2}+\frac{(\log q)^2}{3p}}{p\left(1+\frac{p \log p}{\log q}+\frac{\log p}{2\log q}\right)}\\
&=1-\frac{\log q}{6p^2}\left(\frac{-3p+3p\log q-2(\log q)^2}{\log q+p\log p+\log p/2}\right).\tag{3.6}
\end{align*}
It is easy to verify that
\begin{equation}\tag{3.7}
\frac{-3p+3p\log q-2(\log q)^2}{\log q+p\log p+\log p/2}>\frac{2\log q}{\log p}, \forall\, q\geq 23.
\end{equation}
Combine (3.6) and (3.7), we get
\begin{equation}\tag{3.8}
\frac{G(n_i )}{G(n_{i+1} )}<1-\frac{\log q}{6p^2}\left(\frac{2\log q}{\log p}\right)=1-\frac{(\log q)^2}{3p^2\log p},\quad \forall\, q\geq 23.
\end{equation}
3) Assume $\epsilon_{i+1}=F(q,k)=F(r,j)$, $q\geq 23,r\geq 23$. Then $n_{i+1}=n_i qr$. 
Compare $G(n_i )$ and $G(n_{i+1} )$, we have
\begin{align*}
\frac{G(n_i )}{G(n_{i+1} )}&=\frac{\rho(n_i)\log \log n_{i+1}}{\rho(n_{i+1})\log \log n_i}\\
&=\frac{\log \log n_i+\log \left(1+\frac{\log q+\log r}{\log n_i}\right)}{\log \log n_i}\\
&\quad \times \left(1-\frac{1}{1+q+\cdots+q^k}\right)\left(1-\frac{1}{1+r+\cdots+r^j}\right)\\
&<\left(1+\frac{\log\left(1+\frac{\log q}{p}\right)}{\log p}\right)\left(1+\frac{\log\left(1+\frac{\log r}{p}\right)}{\log p}\right)\\
&\quad \times\left(1-\frac{1}{1+q+\cdots+q^k}\right)\left(1-\frac{1}{1+r+\cdots+r^j}\right).\tag{3.9}
\end{align*}
By Lemma 1 (L1.2'), we have
\begin{equation}\tag{3.10}
q+\cdots+q^k<\frac{1}{\log q}\left(p\log p+\frac{\log p}{2}\right),\,r+\cdots+r^j<\frac{1}{\log r}\left(p\log p+\frac{\log p}{2}\right).
\end{equation}
Then we can proceed with q and r separately as in 2) to prove (3.2).  
\end{proof}

We will prove Lemmas 5-7, then use them to prove an upper bound for $G(n_{i+1} )/G(n_i )$ in Theorem 4. 

\noindent {\bfseries Lemma 5.}
\textit{Define
\begin{equation}\tag{L5.1}
f(x):=\frac{1}{\sqrt{2x}}\sum_{k=3}^{K(x)}(kx)^{1/k},\, x>2.667,
\end{equation}
where $K(x)$ is implicitly defined as the largest integer K satisfying
\begin{equation}\tag{L5.2}
\frac{2^K}{K}\leq x.
\end{equation}
Then\\
1) $f(x)$ is a piece-wise differentiable function with discontinuous points at $x=\frac{2^K}{K}$ for each integer $K\geq 3$.\\
2) $f(x)$ decreases at differentiable points.\\
3)  $f(x)$  has local maximums at discontinuous points $x=\frac{2^K}{K}$.  $f\left(\frac{2^K}{K}\right)>f\left(\frac{2^{K+1}}{K+1}\right)$, for $K\geq 7$. \\
4) In particular, 
\begin{equation}\tag{L5.3}
f(x)<0.10924,\,\forall\, x\geq \frac{2^{31}}{31}=6.93\times 10^7.
\end{equation}
}
\begin{proof}
1) and 2) are simple.
\begin{align*}
f'(x)&=\frac{1}{\sqrt{2}}\left(\sum_{k=3}^{K(x)}k^{\frac{1}{k}}x^{\frac{1}{k}-\frac{1}{2}}\right)'\\
&=\frac{1}{\sqrt{2}}\sum_{k=3}^{K(x)}\left(\frac{1}{k}-\frac{1}{2}\right) k^{\frac{1}{k}}x^{\frac{1}{k}-\frac{1}{2}-1}<0, \,\forall \,x>2.667,\,x\neq\frac{2^K}{K}.\tag{L5.4}
\end{align*}
So f(x) decreases at all differentiable points.\\
3) Because $f(x)$ adds an extra summand 2 at point $x=\frac{2^K}{K}$, it is discontinuous there. To show $f(x)$ decreases from one discontinuous point to next, let $x_s=2^{K(x_s )} /K(x_s )$, i.e. $(K(x_s ) x_s )^{1/K(x_s ) }=2$. Then the next discontinuous point is $x_t:=2^{K(x_t )} /K(x_t )$  where $K(x_t )=K(x_s )+1$. So we have
\begin{equation}\tag{L5.5}
x_s=\frac{2^{K(x_s)}}{K(x_s)}=\frac{2^{K(x_s)+1}}{2K(x_s)}=\frac{2^{K(x_t)}}{2(K(x_t)-1)}=\frac{K(x_t)}{2(K(x_t)-1)}x_t.
\end{equation}
Now we want to show $f(x_s )>f(x_t )$.
\begin{align*}
f(x_s)-f(x_t)&=\frac{1}{\sqrt{2x_s}}\sum_{k=3}^{K(x_s)}(kx_s)^{\frac{1}{k}}-\frac{1}{\sqrt{2x_t}}\sum_{k=3}^{K(x_t)}(kx_t)^{\frac{1}{k}}\\
&=\frac{1}{\sqrt{2}}\sum_{k=3}^{K(x_s)}k^{\frac{1}{k}}\left(x_s^{\frac{1}{k}-\frac{1}{2}}-x_t^{\frac{1}{k}-\frac{1}{2}}\right)-\frac{1}{\sqrt{2x_t}}(K(x_t)x_t)^{\frac{1}{K(x_t)}}\\
&\geq \frac{1}{\sqrt{2}}3^{\frac{1}{3}}\left(x_s^{-\frac{1}{6}}-x_t^{-\frac{1}{6}}\right)-\frac{2}{\sqrt{2x_t}}\\
&= \frac{1.44}{\sqrt{2}}\left(\left(\frac{K(x_t)}{2(K(x_t)-1)}x_t\right)^{-\frac{1}{6}}-x_t^{-\frac{1}{6}}\right)-\frac{2}{\sqrt{2x_t}}\\
&= \frac{1.44}{\sqrt{2}}x_t^{-\frac{1}{6}}\left(\left(\frac{2(K(x_t)-1)}{K(x_t)}\right)^{\frac{1}{6}}-1\right)-\frac{2}{\sqrt{2x_t}}.\tag{L5.6}
\end{align*}
When $K(x_t )\geq 8$, 
\begin{equation}\tag{L5.7}
\left(\frac{2(K(x_t)-1)}{K(x_t)}\right)^{\frac{1}{6}}\geq \left(\frac{14}{8}\right)^{\frac{1}{6}}=1.0977.
\end{equation}
\begin{align*}
f(x_s)-f(x_t)&\geq\frac{0.14077}{\sqrt{2}}x_t^{-\frac{1}{6}}-\frac{2}{\sqrt{2x_t}}\\
&=\frac{2}{\sqrt{2x_t}}\left(0.07038\times x_t^{\frac{1}{3}}-1\right)>0, \,\forall\,x_t\geq 2868.\tag{L5.8}
\end{align*}
For $x_t<2868$ and $K\geq 7$, $f(2^K/K)>f(2^{K+1}/(K+1))$ can be directly calculated:
\begin{longtable}{| r | r | r |}
\hline
K & $x=2^K/K$ & $f(x)$\\
\hline
3	&2.67	&0.87\\
\hline
4	&4.00	&1.52\\
\hline
5	&6.40	&1.94\\
\hline
6	&10.67	&2.15\\
\hline
\textbf{7}	&\textbf{18.29}	&\textbf{2.21}\\
\hline
8	&32.00	&2.16\\
\hline
9	&56.89	&2.03\\
\hline
10	&102.40	&1.86\\
\hline
11	&186.18	&1.67\\
\hline
12	&341.33	&1.48\\
\hline
13	&630.15	&1.30\\
\hline
14	&1170.29	&1.12\\
\hline
15	&2184.53	&0.97\\
\hline
16	&4096.00	&0.83\\
\hline
\end{longtable}
4) Direct calculation shows $f \left(\frac{2^{31}}{31}\right)=0.10923475$.   
\end{proof}

\noindent {\bfseries Lemma 6.}
\textit{    Let $\theta(x)$ and $\psi(x)$ be Chebyshev functions. Define
\begin{equation}\tag{L6.1}
\psi_0 (x):=\sum_{k=1}^K\theta((kx)^{1/k}) ,
\end{equation}
where K is the largest integer k such that $(kx)^{1/k}\geq 2$. Then
\begin{equation}\tag{L6.2}
\psi_0 (x)<x\left(1+\frac{0.06323}{(\log x)^2}\right),\,\forall\, x>10^8.
\end{equation}
}
\begin{proof}
By [PT 2018] Theorem 1, 
\begin{equation}\tag{L6.3}
\theta(x)<x,\quad \forall\, 0<x<1.39\times 10^{17}.
\end{equation}
Setting $k=2,\eta_2=0.01$ in Theorem 4.2 of [Dusart 2018], we have,
\begin{equation}\tag{L6.4}
\vert \theta(x)-x\vert<\frac{0.01x}{(\log x)^2} ,\quad \forall\, x\geq 7\,713\,133\, 853.
\end{equation}

Combine (L6.3) and (L6.4), we get
\begin{equation}\tag{L6.5}
\theta(x)-x<\frac{0.01x}{\log (1.39\times10^{17})^2} =6.418\times 10^{-6}x,\quad \forall \,x>0.
\end{equation}
By (L6.5) and Lemma 5, we have
\begin{align*}
\psi_0 (x)&=\sum_{k=1}^K\theta((kx)^{1/k} ) \\
&<\theta(x)+(1+6.418\times 10^{-6} )\left((2x)^{1/2}+\cdots+(Kx)^{1/K} \right)\\
&<\theta(x)+(1+6.418\times 10^{-6} ) (2x)^{1/2} (1+0.109235)\\
&<x+\frac{0.01x}{(\log x)^2} +1.5687x^{1/2}\\                                                
&=x+\frac{x}{(\log x)^2} \left(0.01+\frac{1.5687(\log x )^2}{x^{1/2}}                      \right)\\
&<x+\frac{0.06323x}{(\log x )^2},\quad     \forall\, x>10^8.\tag{L6.6}
\end{align*}
\end{proof}

\noindent {\bfseries Lemma 7.}
\textit{Let n be a CA number and $p=P(n)$. Then
\begin{equation}\tag{L7.1}
\log n<p\left(1+\frac{0.06323}{(\log p )^2} \right),\quad \forall\, p>10^8.
\end{equation}
}
\begin{proof}
Let $x_k$ be defined by (4). By method of Theorem 4 of [Wu 2019], we have
\begin{align*}
\log N&=\theta(p)+\theta(x_2 )+\cdots+\theta(x_K )\\
&<\theta(p)+\theta((2p)^{1/2} )+\cdots+\theta((Kp)^{1/K} )=\psi_0 (p),\tag{L7.2}
\end{align*}
where K is the largest integer k such that $(kx)^{1/k}\geq 2$ and $\psi_0$ is defined as in Lemma 6. By Lemma 6, we have
\begin{equation}\tag{L7.3}
\log N<p\left(1+\frac{0.06323}{(\log p)^2}\right),\quad \forall\, p>10^8.
\end{equation}
\end{proof}

\noindent {\bfseries Theorem 4.}
\textit{Let $n_i$ be CA3. Let p be the prime succeeding $P(n_i )$ and $p>10^8$. Then $p<\log n_i$, Assume $n_{i+1}=n_iq$ and $\epsilon_{i+1}=F(q,k)$ for some prime q and integer k. Then 
\begin{equation}\tag{4.1}
G(n_{i+1} )<G(n_i ) \exp\left(\frac{0.12646 \log q}{p(\log p)^3}\right).
\end{equation}
}
\begin{proof}
Write $\epsilon:=\epsilon_{i+1}$. Define $g(t):=t^\epsilon/\log \log t$  with minimum at $t_0$. Then we have
\begin{align*}
\frac{G(n_{i+1})}{G(n_i )}&=\frac{g(n_{i+1} )}{g(n_i )}=q^\epsilon \frac{\log \log n_i}{\log \log n_{i+1}}\\  
&=\exp(\epsilon \log q+\log \log \log n_i -\log \log \log n_{i+1} )\\
&<\exp \left(\frac{\log q}{\log t_0 \log \log t_0}-\frac{\log \log n_{i+1}-\log \log n_i}{\log \log n_{i+1}} \right)\\
&<\exp \left(\frac{\log q}{\log t_0 \log \log t_0}-\frac{\log n_{i+1}- \log n_i}{\log n_{i+1} \log \log n_{i+1}} \right)\\
&=\exp \left(\frac{\log q}{\log t_0 \log \log t_0}-\frac{\log q}{\log n_{i+1} \log \log n_{i+1}} \right).\tag{4.2}
\end{align*}
By Lemma 7, 
\begin{equation}\tag{4.3}
\log n_{i+1}<\psi_0 (p)<cp,
\end{equation}
$c:=1+\frac{0.06323}{(\log p )^2}$, for $p>10^8$. Hence
\begin{align*}
\frac{c^2-1}{c^2}= \frac{1+\frac{0.12646}{(\log p)^2}+\frac{0.004}{(\log p)^4}-1}{1+\frac{0.12646}{(\log p)^2}+\frac{0.004}{(\log p)^4}}\\
=\frac{0.12646}{(\log p)^2}\times \frac{1+\frac{0.031615}{(\log p)^2}}{1+\frac{0.12646}{(\log p)^2}+\frac{0.004}{(\log p)^4}}<\frac{0.12646}{(\log p)^2}.\tag{4.4}
\end{align*}
By Lemma 2, $p<\log t_0$. So (4.2) means
\begin{align*}
\frac{G(n_{i+1} )}{G(n_i )}&<\exp\left(\frac{\log q}{p \log p}-\frac{\log q}{c^2 p \log p} \right)\\
&=\exp\left(\frac{(c^2-1) \log q}{c^2 p \log p}\right)<\exp\left(\frac{0.12646 \log q}{p(\log p )^3}\right), \quad \forall  \,p>10^8.\tag{4.5}
\end{align*}

\end{proof}

\begin{center}
{\bfseries \large References}
\end{center}

\noindent {[}Briggs 2006{]} K. Briggs. \textit{Abundant numbers and the Riemann hypothesis}. Experiment. Math., 15(2):251–256, 2006.\\
{[}Broughan 2017{]} K. Broughan, \textit{Equivalents of the Riemann Hypothesis} Vol 1. Cambridge Univ. Press. (2017)\\
{[}CLMS 2007{]} Y. -J. Choie, N. Lichiardopol, P. Moree, and P. Sol\'{e}. \textit{On Robin’s criterion for the Riemann hypothesis}. J. Th\'{e}or. Nombres Bordeaux, 19(2):357–372, 2007.\\
{[}CNS 2012{]} G. Caveney, J.-L. Nicolas, and J. Sondow, \textit{On SA, CA, and GA numbers}, Ramanujan J. 29 (2012), 359–384.\\
{[}Dusart 1998{]} P. Dusart \textit{Sharper bounds for $\psi$, $\theta$, $\pi$, $p_k$}, Rapport de recherche $n\,1998-06$, Laboratoire d’Arithm\'{e}tique de Calcul formel et d’Optimisation\\
{[}Dusart 2018{]} P. Dusart. \textit{Explicit estimates of some functions over primes}. Ramanujan J., 45(1):227–251, 2018.\\
{[}EN 1975{]} P. Erdös and J.-L Nicolas, \textit{Repartition des nombres superabundants}. Bulletin de la S. M., tome 103 (1975), p. 65-90 \\
{[}Morrill;Platt 2018{]} T. Morrill, D. Platt. \textit{Robin’s inequality for 25-free integers and obstacles to analytic improvement} \\
https://arxiv.org/abs/1809.10813\\
{[}NY 2014{]} S. Nazardonyavi and S. Yakubovich. \textit{Extremely Abundant Numbers and the Riemann Hypothesis}. Journal of Integer Sequences, Vol. 17 (2014), Article 14.2.8\\
{[}PT 2016{]} D.J. Platt and Tim Trudgian. \textit{On the first sign change of $\theta(x) - x$}. Math. Comp. 85 (2016), 1539-1547\\
{[}Robin 1984{]} G. Robin. \textit{Grandes valeurs de la fonction somme des diviseurs et hypoth\'{e}se de Riemann}. Journal de mathématiques pures et appliqu\'{e}es. (9), 63(2):187–213, 1984.\\
{[}Wu 2019{]} X. Wu. \textit{Properties of counterexample of Robin hypothesis}.\\
https://arxiv.org/abs/1901.09832

\end{document}